\documentclass[12pt,reqno]{amsart}

\usepackage{amsmath,amsthm,amsfonts,amscd,amssymb,euscript,enumerate}
\numberwithin{equation}{section} 
\newcommand{\dbar}{\ensuremath{\bar \partial}}

\newcommand{\C}{\ensuremath{{\mathbb C}}}
\newcommand{\pj}{\ensuremath{{\mathbb P}}}
\newcommand{\R}{\ensuremath{{\mathbb R}}}
\newcommand{\N}{\ensuremath{{\mathbb N}}}

\newcommand{\smooth}{\ensuremath{C^{\infty}}}

\newcommand{\I}{\ensuremath{\mathcal I}}

\newcommand{\oka}{\ensuremath{\mathcal O}}

\newcommand{\cv}{\ensuremath{\mathcal C}}

\newtheorem{Lemma}{Lemma}
\newtheorem{Theorem}[Lemma]{Theorem}

\newtheorem{Proposition}[Lemma]{Proposition}
\newtheorem{Corollary}[Lemma]{Corollary}
\newtheorem*{Theorem*}{Theorem}
\newtheorem*{thm1.1}{Theorem 1.1}
\newtheorem*{thm1.2}{Theorem 1.2}
\newtheorem*{thm2.9'}{Theorem 2.9'}
\newtheorem*{prop4.1}{Proposition 4.1}
\newtheorem*{prop4.2}{Proposition 4.2}
\newtheorem*{prop4.4}{Proposition 4.4}
\newtheorem*{prop4.5}{Proposition 4.5}
\theoremstyle{remark}
\newtheorem{Remark}[Lemma]{Remark}

\newtheorem*{Remark*}{Remark}
\theoremstyle{definition}
\newtheorem{Definition}[Lemma]{Definition}
\newtheorem*{Definition*}{Definition}
\newtheorem{Example}[Lemma]{Example}

\numberwithin{Lemma}{section}

\numberwithin{equation}{section}

\makeatletter
\@namedef{subjclassname@2020}{%
  \textup{2020} Mathematics Subject Classification}
\makeatother

\begin{document}
\title[Notions of Higher Type]{Notions of Higher Type}

\author{Andreea C. Nicoara}

\address{School of Mathematics, Trinity College Dublin, Dublin 2, Ireland}

\email{anicoara@maths.tcd.ie}



\begin{abstract}
Notions of finite type play an important role in several complex variables. The most standard notion is D'Angelo type, which measures the order of contact of holomorphic curves with the boundary of a domain in $\C^n.$ For the $\dbar$-Neumann problem, however, the order of contact of the boundary of the domain with $q$-dimensional complex varieties controls its behavior on $(p,q)$ forms. There are two different ways of measuring this order of contact, one due to John D'Angelo and another due to David Catlin. We survey the known results about the D'Angelo and Catlin $q$-types, their relationship, and other notions that complete the picture.
\end{abstract}

\dedicatory{Dedicated to Vasile Brinzanescu on his 80th birthday}

\maketitle

\tableofcontents

\section{Introduction}
\label{intro}

Holomorphic functions by their very definition sit in the kernel of the $\dbar$ operator making the $\dbar$-Neumann problem the most important partial differential equation for the field of several complex variables. Joseph J. Kohn solved the $\dbar$-Neumann problem in \cite{dbneumann1} and \cite{dbneumann2} on strongly pseudoconvex domains in $\C^n$ and proved that for a strongly pseudoconvex domain the solution only gains half of the order of the operator in the Sobolev norm. Thus, at best the $\dbar$-Neumann problem is subelliptic. Kohn's result launched a quest to characterize when the $\dbar$-Neumann problem is subelliptic on pseudoconvex domains, which are the most natural objects in complex analysis of several variables as these are exactly the domains of holomorphy. \cite{hormanderbook} and \cite{krantzbook} are excellent sources on the background material pertaining to domains of holomorphy, pseudoconvexity, and the Levi problem. In analogy to H\"{o}rmander's work in \cite{hormandersquares} as well as of Nirenberg and Tr\`{e}ves in \cite{nirenbergtreves1} and \cite{nirenbergtreves2} and of Tr\`{e}ves in \cite{treves70}, in 1972  in \cite{kohn72} Joseph J. Kohn introduced a finite commutator type condition on a pseudoconvex domain in $\C^2$ and proved that this condition was exactly the right one in order to ensure subellipticity of the $\dbar$-Neumann problem for $(0,1)$ forms on a pseudoconvex domain in $\C^2.$ Dimension 2 is special, however, because the Levi form is a function so the theory is significantly easier. In 1977 in \cite{bloomgrahamjdg}, Thomas Bloom and Ian Graham generalized Kohn's commutator type condition and proved that for a real hypersurface in $\C^n$ it is equal to the order of contact of $(n-1)$-dimensional manifolds with the hypersurface regardless of whether pseudoconvexity holds or not. This commutator type condition is now called the Bloom-Graham type. 

Fifteen years after he established the subellipticity of the $\dbar$-Neumann problem on strongly pseudoconvex domains, Joseph J. Kohn proved that the singularities of $\dbar$ propagate along complex varieties in \cite{kohnacta} and thus was able to tease out the correct geometric condition on a pseudoconvex domain in order to ensure behavior comparable to that of strongly pseudoconvex domains for the $\dbar$-Neumann problem on $(p,q)$ forms: finite order of contact of complex varieties of dimension $q$ with the boundary of the domain. Using a fair amount of algebraic geometry and a result of Klas Diederich and John Erik Forn\ae ss in \cite{df}, Kohn showed that for a pseudoconvex domain with real analytic boundary, the subellipticity of the $\dbar$-Neumann problem on $(p,q)$ forms is equivalent to the condition that complex varieties of dimension $q$ have finite order of contact with the boundary of the domain. One of the examples that Kohn uses to illustrate the subtleties of type in \cite{kohnacta} is given by $$r=\text{Re} \, {z_1} + |z_2^2-z_3^3|^2,$$ which defines a pseudoconvex domain with real analytic boundary. The order of contact of $2$-dimensional complex manifolds with the boundary of this domain at the origin is 4, which equals the Bloom-Graham type of this domain at the origin. The order of contact of $1$-dimensional complex manifolds with the boundary of this domain at the origin is 6. The $\dbar$-Neumann problem for $(0,1)$ forms on this domain is not subelliptic at the origin, however, because the complex curve $z_2^2=z_3^3$ is sitting in the boundary of this domain so the order of contact of the boundary of the domain with $1$-dimensional complex varieties, namely complex curves, is infinite. The commutator type is finite at the origin because the complex curve $z_2^2=z_3^3$ has a singularity there, which commutators of vector fields cannot pick up.

An issue left open by Kohn's work in \cite{kohnacta} was whether the equivalence between the subellipticity of the $\dbar$-Neumann problem on $(p,q)$ forms with the finite order of contact of $q$-dimensional complex varieties with the boundary of the domain also held in even more natural setting of a pseudoconvex domain with $\smooth$ boundary in $\C^n.$ Kohn's approach in the real analytic case was to construct subelliptic multipliers via an algebraic geometric procedure involving taking real radicals and forming determinants with the aim of proving the existence of a non-zero subelliptic multiplier that would establish the needed subelliptic estimate. Such procedure runs into a significant issue in the $\smooth$ case, namely that no Nullstellensatz result is known with respect to the real radical or any other notion of radical. As a result, this problem is still open and is called the Kohn Conjecture. 

The notions of type explored by Thomas Bloom and Ian Graham, while not so useful for characterizing subellipticity of the $\dbar$-Neumann problem, have had an interesting development of their own. In analogy to a proof given by Kohn in \cite{kohnacta} for subellipticity of the $\dbar$-Neumann problem on $(p,n-1)$ forms, Bloom introduced in \cite{bloomsns} a commutator type condition on partial traces of the Levi form now called the Levi form type and showed that for a pseudoconvex domain the Levi form type on the entire trace of the Levi form equals the Bloom-Graham type of the boundary of the domain and also the order of contact of the boundary of the domain with $(n-1)$-dimensional complex manifolds. Thus, in \cite{bloom81} Bloom conjectured that commutator type defined using $q$-dimensional subbundles of the $(1,0)$ tangent space of a pseudoconvex hypersurface in $\C^n$ equals the order of contact of the hypersurface with $q$-dimensional complex manifolds and also the Levi form type at level q. This statement is called the Bloom conjecture. It was proven for pseudoconvex hypersurfaces in $\C^3$ in 2021 by Xiaojun Huang and Wanke Yin in \cite{huangyin}, whereas for $\C^n$ with $n \geq 4$ only partial results are known.

Meanwhile, in 1982 in \cite{opendangelo} John D'Angelo formalized the notion of order of contact of varieties of dimension $q$ with a hypersurface in $\C^n$ such as the boundary of a domain. This notion $\Delta_q$ is now called the D'Angelo $q$-type, and specifically for $q=1,$ it is the most standard notion of finite type in the several complex variables literature. We will devote part of this article to explaining its properties in detail.

In the mid 1980's after D'Angelo's work, David Catlin set out to prove that the subellipticity of the $\dbar$-Neumann problem  on $(p,q)$ forms on a pseudoconvex domain with $\smooth$ boundary in $\C^n$ is equivalent to finite $q$-type in a manner that used the geometry of the domain and more standard analytic methods in several complex variables but avoided the algebraic geometric arguments employed by Kohn in \cite{kohnacta} for the real analytic case. Catlin first proved the necessity of the finite $q$-type condition in 1983 in \cite{catlinnec}. The sufficiency of the condition took considerably more work. First, in 1984 in \cite{catlinbdry} Catlin sought to define a notion that captured the vanishing order of the defining function of the domain in various directions, which he called the multitype. He gave a rather abstract definition of the multitype and then defined a more concrete notion using lists of vector fields that he called the commutator multitype. For a pseudoconvex domain in $\C^n,$ Catlin gave a very intricate proof that the multitype equals the commutator multitype that has yet to be fully understood. In private communication with the author, Catlin credited drawing inspiration for the commutator multitype from the work of Bloom and Graham in \cite{bloomgrahamjdg} and \cite{bloomgrahaminv} though Kohn's arguments in section 5 of \cite{kohnacta} also clearly had an impact on his thinking. In the main theorem of \cite{catlinbdry}, Catlin states that the $\nu^{th}$ entry of the multitype for $2 \leq \nu \leq n$ is controlled by the D'Angelo $(n-\nu+1)$-type $\Delta_{n-\nu+1},$ but from the proof it is clear that it is rather the regular $(n-\nu+1)$-type, namely the order of contact of $(n-\nu+1)$-dimensional manifolds with the boundary of the domain rather than varieties, that is the upper bound for the $\nu^{th}$ entry of the multitype. This fact could perhaps be used to produce another approach for proving or disproving the Bloom Conjecture when $n \geq 4.$ Finally, in 1987 in \cite{catlinsubell}, Catlin tackled the sufficiency of the finite $q$-type condition for the subellipticity of the $\dbar$-Neumann problem  on $(p,q)$ forms on a pseudoconvex domain with $\smooth$ boundary in $\C^n.$ Rather than using D'Angelo's notion of $q$-type $\Delta_q,$ Catlin defined another notion in a far more complicated but more generic fashion, $D_q,$ which is now called the Catlin $q$-type. The two notions, $\Delta_q$ and $D_q$ trivially equal each other when $q=1$ and were presumed to be equal to each other for higher $q$ as well when the author and Vasile Brinzanescu started thinking about their relationship in 2012 during Vasile Brinzanescu's semester long visit to the University of Pennsylvania; see the survey \cite{kohndangelo} published by D'Angelo and Kohn in 1999.

The 1982 paper \cite{opendangelo} where John D'Angelo defines $\Delta_q$ mostly deals with $\Delta_1.$ The higher type only occupies a section of the paper at the very end. Likewise D'Angelo's 1993 monograph \cite{dangelo} does not devote more than a section to $\Delta_q$ for $q>1.$ Such was the state of affairs in 2012. Vasile Brinzanescu and the author published a first paper \cite{bazilandreea} on $q$-types in 2015 deriving effective relationships between the Catlin and the D'Angelo $q$-types without being able to prove equality of the two notions. In 2019 Martino Fassina, a PhD student of John D'Angelo at the time, published an example in \cite{fassina} that showed not only that $\Delta_q$ and $D_q$ could be different but also that examples of ideals of holomorphic functions could be constructed so that these two notions would differ by any chosen positive integer. In 2023 Vasile Brinzanescu and the author gave a much simpler to understand and compute restatement of the Catlin $q$-type and fully elucidated the relationship between the Catlin and D'Angelo $q$-types in \cite{bazilandreea2}. A further work is in preparation that will prove that both the Catlin and the D'Angelo $q$-types are biholomorphic invariants, giving an affirmative answer to a question posed by Dmitri Zaitsev to John D'Angelo and to the author after the publication of \cite{bazilandreea}. 

The paper is organized as follows: Section~\ref{findangelo} is devoted to the D'Angelo $q$-type, section~\ref{fincatlin} to the Catlin $q$-type, and section~\ref{relate} to their relationship. 

\section{D'Angelo $q$-type}
\label{findangelo}

The aim of this section is to give a brief overview of the properties of the D'Angelo $q$-type. For full details, the reader is directed to D'Angelo's 1982 paper \cite{opendangelo} where he defines his $q$-type notion as well as to his 1993 monograph \cite{dangelo}.

Holomorphic curves have parametrizations so measuring the order of contact of a hypersurface in $\C^n$ with holomorphic curves is easy to do as we can pull back the defining function of the hypersurface by the parametrization of the curve and measure its order of vanishing. Holomorphic curves do not have unique parametrizations, however, so we need to normalize by the vanishing order of the holomorphic curve. This idea underpins D'Angelo's definition of $\Delta_1,$ his $1$-type. More precisely, let $r$ be a defining function for a real hypersurface $M$ in $\C^n,$ let $x_0 \in M,$ and let $\cv=\cv(n,p)$ be the set of all germs of holomorphic curves $$\varphi: (U,0) \rightarrow (\C^n, x_0),$$ where $U$ is some neighborhood of the origin in $\C^1$ and $\varphi(0)=x_0.$ Such a curve $$\varphi(t)=(\varphi_1(t), \dots, \varphi_n(t))$$ satisfies that $\varphi_j(t)$ is holomorphic for every $j$ with $1 \leq j \leq n$ and every $t \in U.$ For each component  $\varphi_j(t),$ let ${\text ord}_0 \, \varphi_j$ be the order of the first non-vanishing derivative of $\varphi_j$ at the origin, and define ${\text ord}_0 \, \varphi = \min_{1 \leq j \leq n} \, {\text ord}_0 \, \varphi_j.$ By convention, if $\varphi_j(0) \equiv 0,$ then its order of vanishing at the origin is $+\infty,$ whereas if $\varphi_j(0) \neq 0,$ then its order of vanishing at the origin is $0.$

\medskip
\newtheorem{firstft}{Definition}[section]
\begin{firstft}
Let $M$ be a real hypersurface in $\C^n,$ and let $r$ be a defining function for $M.$ The D'Angelo $1$-type at $x_0 \in M$ is given by $$\Delta_1 (M, x_0) = \sup_{\varphi \in \cv(n,x_0)} \frac{ {\text ord}_0 \, \varphi^* r}{{\text ord}_0 \, \varphi},$$ where $\varphi^* r,$ the pullback of $r$ to $\varphi$ and $ {\text ord}_0 \, \varphi^* r$ is its order of vanishing at the origin.
\end{firstft}

\medskip\noindent {\bf Remark:}

\noindent Kohn's commutator type, the Bloom-Graham $q$ commutator type, and the Levi form type at level $q$ are all positive integers and upper semi-continuous by their very definitions. $\Delta_1$ is positive but does not have to be an integer due to the normalization by the order of the curve and is not upper semi-continuous. We shall now provide examples to illustrate both of these statements.

\begin{Example}
We modify an example of D'Angelo from p.73 of  \cite{dangelo} to illustrate that $\Delta_1$ is not necessarily an integer. Let the hypersurface $M$ have defining function given by $$r= 2 \text{Re} \, {z_3} + |z_1^3-z_2^5|^2 +|z_1^2 z_2^3|^2.$$ There are only three relevant curves as far as computing $\Delta_1.$ The first curve is $\varphi(t)=(t^5,t^3,0),$ which gives ${\text ord}_0 \, \varphi^* r=38,$ whereas ${\text ord}_0 \, \varphi=3.$ The second curve is $\varphi(t)=(t,0,0),$ which gives ${\text ord}_0 \, \varphi^* r=6$ and ${\text ord}_0 \, \varphi=1.$ Finally, the third curve is $\varphi(t)=(0,t,0),$ which gives ${\text ord}_0 \, \varphi^* r=10$ and ${\text ord}_0 \, \varphi=1.$ Clearly, the first curve achieves the supremum and gives $\Delta_1 (M,0)=\frac{38}{3},$ which is not an integer.
\end{Example}

\begin{Example}
\label{dangeloex2}
The following example due to D'Angelo illustrates the failure of $\Delta_1$ to be upper semi-continuous: Let the hypersurface $M$ have defining function given by $$r= 2 \text{Re} \, {z_3} + |z_1^2-z_2 z_3|^2 +|z_2^2|^2,$$ then $\Delta_1(M, 0)=4$ because any curve that has a non-zero third component will pick up the term $2 \text{Re} \, {z_3}$ and not return the supremum whereas any curve with either a non-zero second component or a non-zero third component will give $4$ in the quotient $\displaystyle  \frac{ {\text ord}_0 \, \varphi^* r}{{\text ord}_0 \, \varphi}.$ Now, let us look at points of the form $x=(0,0,ia),$ where $a \in \R.$ We have $\Delta_1(M, x) = 8$ because the supremum is achieved by the holomorphic curve $\displaystyle \varphi(t)=\left(t, \frac{t^2}{i \epsilon}, i \epsilon \right),$ which annihilates $ 2 \text{Re} \, {z_3} + |z_1^2-z_2 z_3|^2.$  For more details, see p.136 of \cite{dangelo} as well as \cite{nonusc}. The next theorem from \cite{opendangelo} shows that this example exhibits the maximum possible jump in $1$-type.
\end{Example}

\begin{Theorem}[Theorem 5.7 p.632 \cite{opendangelo}]
Let a hypersurface $M$ be pseudoconvex at a point $x_0,$ and let the Levi form have rank $p$ there. For every $x \in M$ sufficiently close to $x_0,$ $$\Delta_1 (M,x) \leq \frac{\Delta_1(M,x_0)^{n-1-p}}{2^{n-2-p}}.$$
\end{Theorem}

\noindent In the previous example, $n=3$ and $p=0$ at $x_0=0.$ For every $x=(0,0,ia),$ $$ \Delta_1(M, x) = 8 =\frac{4^{3-1}}{2^{3-2}}.$$

\medskip\noindent To define $\Delta_q$ for $q >1,$ we observe that $q$-dimensional complex varieties do not necessarily have parametrizations so it is necessary to reduce this situation to some type of computation involving $\Delta_1.$ To obtain holomorphic curves starting with complex varieties of dimension $q,$ we can cut the ambient space $\C^n$ using a set of $q-1$ linearly independent linear forms $w_1, \dots, w_{q-1}.$ The zero locus of these forms gives a linear embedding of $\C^{n-q+1}$ into $\C^n,$ which we can denote by $\phi : \C^{n-q+1} \rightarrow \C^n.$ Generically, namely for generic sets of linearly independent linear forms $w_1, \dots, w_{q-1},$ the pullback $\phi^* M$ is a hypersurface in $\C^{n-q+1}.$ D'Angelo chose to define $\Delta_q$ for $q>1$ in \cite{opendangelo} as follows:

\medskip
\newtheorem{qft}[firstft]{Definition}
\begin{qft}
\label{qfinitetype}
Let $M$ be a real hypersurface in $\C^n,$ and let $r$ be a defining function for $M.$ The D'Angelo $q$-type at $x_0 \in M$ is given by $$\Delta_q (M, x_0) =\inf_\phi \sup_{\varphi \in \cv(n-q+1,x_0)} \frac{ {\text ord}_0 \, \varphi^* \phi^*r}{{\text ord}_0 \, \varphi}=\inf_\phi \Delta_1 (\phi^*r, x_0),$$ where $\phi : \C^{n-q+1} \rightarrow \C^n$ is any linear embedding of $\C^{n-q+1}$ into $\C^n$ and we have identified $x_0$ with $\phi^{-1} (x_0).$ 
\end{qft}

\medskip\noindent {\bf Remark:}

\noindent It is essential to note here that D'Angelo made the choice of taking the infimum over all possible linear embeddings $\phi,$ and in algebraic geometry the infimum is not a very robust notion. Catlin opted to define his own notion of $q$-type $D_q$ in \cite{catlinsubell} precisely in order to fix this particular issue, even though his definition ended up being quite complicated both to compute and to understand.

\smallskip\noindent The next theorem assembles together the most important properties of $\Delta_q:$

\medskip
\newtheorem{propdeltaq}[firstft]{Theorem}
\begin{propdeltaq}
Let $M$ be a smooth real hypersurface in $\C^n.$
\begin{enumerate}
\item[(i)] $\Delta_q (M, x_0)$ is well-defined, i.e. independent of the defining function $r$ chosen for $M.$
\item[(ii)] $\Delta_q (M, x_0)$ is not necessarily upper semi-continuous.
\item[(iii)]  Let $\Delta_q (M, x_0)$ be finite at some $x_0 \in M,$ then there exists a neighborhood $V$ of $x_0$ on which $$\Delta_q (M, x) \leq 2 (\Delta_q(M, x_0))^{n-q}.$$
\item[(iv)] The function $\Delta_q(M, x_0)$ is finite determined. In other words, if $\Delta_q(M, x_0)$ is finite, then there exists an integer $k$ such that $\Delta_q(M, x_0)=\Delta_q(M',x_0)$ for $M'$ a hypersurface defined by any $r'$ that has the same $k$-jet at $x_0$ as the defining function $r$ of $M.$
\end{enumerate} \label{propdeltaqthm}
\end{propdeltaq}

\medskip\noindent {\bf Remarks:}

\noindent (1) Part (iii) is D'Angelo's Theorem 6.2 from p.634 of \cite{opendangelo} and shows the set of points of finite $q$-type is open. Pseudoconvexity is not required for this conclusion.

\noindent (2) Part (iv) is Proposition 14 from p.88 of \cite{dangeloftc}. When we examine the proof of this result, we see that if $t = \Delta_q (M, x_0) < \infty,$ then we can take $k = \lceil t \rceil,$ the ceiling of $t,$ which is by definition the least integer greater than or equal to $t.$ We direct the reader to p.628 of \cite{opendangelo} for an examination on a concrete example of how the type of the truncation can change when $k<\lceil t \rceil.$

\medskip\noindent After his work in \cite{opendangelo}, D'Angelo realized that it would be helpful to extend the notions of $\Delta_1$ and $\Delta_q$ to ideals rather than dealing only with hypersurfaces in $\C^n.$ The following two definitions come from p.86 of \cite{dangeloftc}:

\medskip
\newtheorem{idtype}[firstft]{Definition}
\begin{idtype}
\label{idealtype}
Let $C^\infty_{x_0}$ be the ring of smooth germs at $x_0 \in \C^n$ and let $\I$ be an ideal in $C^\infty_{x_0}.$ $$\Delta_1(\I,x_0)   = \sup_{\varphi \in \cv(n,x_0)} \:\: \inf_{g \in \I} \:\:\frac{ {\text ord}_0 \, \varphi^* g}{{\text ord}_0 \, \varphi}.$$
\end{idtype}

\medskip
\newtheorem{qftequiv}[firstft]{Definition}
\begin{qftequiv}
\label{qfinitetypeequiv}
Let $M$ be a real hypersurface in $\C^n,$ and let $x_0 \in M.$ The D'Angelo $q$-type at $x_0 \in M$ is given by $$\Delta_q (M, x_0) =\inf_{\{w_1, \dots, w_{q-1} \}} \Delta_1\Big( (\I(M), w_1, \dots, w_{q-1}), x_0 \Big),$$ where $\{w_1, \dots, w_{q-1} \}$ is a non-degenerate set of linear forms, $(\I(M), w_1, \dots, w_{q-1})$ is the ideal in $C^\infty_{x_0}$ generated by $\I(M), w_1, \dots, w_{q-1},$ and the infimum is taken over all such non-degenerate sets $\{w_1, \dots, w_{q-1} \}$ of linear forms.\end{qftequiv}

\medskip\noindent {\bf Remarks:}

\noindent (1) This definition is equivalent to Definition~\ref{qfinitetype}.

\noindent (2) Note that the curves $\varphi \in \cv(n,x_0)$ in the definition of $ \Delta_1\Big( (\I(M), w_1, \dots, w_{q-1}), x_0 \Big)$ satisfy that they sit in the complex linear subspace given by $w_1 = \dots = w_{q-1}=0$ because otherwise the infimum in the definition of $ \Delta_1\Big( (\I(M), w_1, \dots, w_{q-1}), x_0 \Big)$ is achieved on one of the $w_j$'s and then the supremum over $\cv(n,x_0)$ cannot be realized.

\medskip\noindent This same kind of definition makes sense for an ideal of holomorphic germs as well, so D'Angelo also defined the following:

\medskip
\newtheorem{qftideal}[firstft]{Definition}
\begin{qftideal}
\label{qfinitetypeideal}
Let $\oka_{x_0}$ be the local ring of holomorphic germs in $n$ variables at $x_0 \in \C^n.$If $\I$ is an ideal in $\oka_{x_0},$
\begin{equation*}
\begin{split}
\Delta_q(\I, x_0)&=\inf_{\{w_1, \dots, w_{q-1} \}} \Delta_1\Big( (\I, w_1, \dots, w_{q-1}), x_0 \Big)\\&=\inf_{\{w_1, \dots, w_{q-1} \}}\sup_{\varphi \in \cv(n,x_0)} \:\: \inf_{g \in (\I, w_1, \dots, w_{q-1})} \:\:\frac{ {\text ord}_0 \, \varphi^* g}{{\text ord}_0 \, \varphi},
\end{split}
\end{equation*} 
where $\{w_1, \dots, w_{q-1} \}$ is a non-degenerate set of linear forms in $\oka_{x_0},$ $(\I, w_1, \dots, w_{q-1})$ is the ideal in $\oka_{x_0}$ generated by $\I, w_1, \dots, w_{q-1},$ and the infimum is taken over all such non-degenerate sets $\{w_1, \dots, w_{q-1} \}$ of linear forms in $\oka_{x_0}.$
\end{qftideal}

In order to prove the local boundedness of type, which implies that the set of points of finite type is open, D'Angelo used polarization of the truncation of the defining function of the hypersurface $M$ and also introduced his property P. Let us start by describing D'Angelo's polarization argument. If $\Delta_q (M, x_0) = t < \infty,$ let $k = \lceil t \rceil.$ We apply Theorem~\ref{propdeltaqthm} (iv) and remark (2) following it to conclude that $\Delta_q (M, x_0) = \Delta_q (M_k, x_0)$ for $M_k$ the real hypersurface with defining function $r_k,$ which is the polynomial whose $k$-jet at $x_0$ is the same as that of the defining function $r$ of $M.$ Next, we give a holomorphic decomposition for $r_k$ as $$r_k = Re\{h\} + ||f||^2-||g||^2,$$ where $||f||^2= \sum_{j=1}^N |f_j|^2,$ $||g||^2= \sum_{j=1}^N |g_j|^2,$ and the functions $h, f_1, \dots, f_N, g_1, \dots, g_N$ are all holomorphic polynomials in $n$ variables; see Section III of \cite{opendangelo}. We can now state D'Angelo's property P, which first appeared on p.631 of \cite{opendangelo}:

\medskip
\newtheorem{propP}[firstft]{Definition}
\begin{propP}
\label{propertyP}
Let $M$ be a real hypersurface of $\C^n,$ and let $x_0$ be a point of finite type on $M.$ We suppose that $\Delta_1(M, x_0) <k.$ Let $j_{k,x_0} r = r_k = Re\{h\} + ||f||^2-||g||^2$ be a holomorphic decomposition at $x_0$ of the $k$-jet of the defining function $r$ of $M.$ We say that $M$ satisfies property P at $x_0$ if for every holomorphic curve $\varphi \in \cv(n,x_0)$ for which $\varphi^* h$ vanishes, the following two conditions are satisfied:
\begin{enumerate}
\item[(i)] ${\text ord}_0 \, \varphi^* r$ is even, i.e. ${\text ord}_0 \, \varphi^* r = 2a,$ for some $a \in \N;$
\item[(ii)] $\displaystyle \left(\frac{d}{dt}\right)^a \left(\frac{d}{d \bar t}\right)^a \varphi^* r (0) \neq 0.$
\end{enumerate}
\end{propP}

\medskip\noindent {\bf Remarks:}  

\noindent (1) Since $\Delta_q(M, x_0)$ is finitely determined for all $1 \leq q <n,$ this definition is independent of $k$ provided $k$ is large enough. 

\noindent (2) The D'Angelo type is always at least $2,$ whereas the function $h$ has order $1$ at $x_0,$ so no holomorphic curve with $\varphi^* h \not\equiv 0$ can realize the supremum. We can thus exclude such curves from the computation of the type.

\bigskip\noindent In \cite{bazilandreea} Vasile Brinzanescu and the author introduced the $q$ version of D'Angelo's property P, $q$-positivity:

\medskip
\newtheorem{qpropP}[firstft]{Definition}
\begin{qpropP}
\label{qpropertyP}
Let $M$ be a real hypersurface of $\C^n,$ and let $x_0 \in M$ be such that $\Delta_q(M, x_0) <k.$ Let $j_{k,x_0} r = r_k = Re\{h\} + ||f||^2-||g||^2$ be a holomorphic decomposition at $x_0$ of the $k$-jet of the defining function $r$ of $M.$ We say that $M$ is $q$-positive at $x_0$ if for every holomorphic curve $\varphi \in \cv(n,x_0)$ for which $\varphi^* h$ vanishes and such that the image of $\varphi$ locally lies in the zero locus of a non-degenerate set of linear forms $\{w_1, \dots, w_{q-1} \}$ at $x_0,$ the following two conditions are satisfied:
\begin{enumerate}
\item[(i)] ${\text ord}_0 \, \varphi^* r$ is even, i.e. ${\text ord}_0 \, \varphi^* r = 2a,$ for some $a \in \N;$
\item[(ii)] $\displaystyle \left(\frac{d}{dt}\right)^a \left(\frac{d}{d \bar t}\right)^a \varphi^* r (0) \neq 0.$
\end{enumerate}
\end{qpropP}

\bigskip\noindent A pseudoconvex domain of finite D'Angelo $1$-type has property P. The hypersurfaces corresponding to truncations of the defining function at $x_0$ of any order higher than the type also have property P. The equivalent statements are true for $q$-positivity, namely a pseudoconvex domain of finite D'Angelo $q$-type is $q$-positive. All hypersurfaces that correspond to truncations of the defining function of the domain of order higher than the type are $q$-positive; see \cite{bazilandreea}.

\section{Catlin $q$-type}
\label{fincatlin}

In the definition of the D'Angelo $q$-type, the ambient space $\C^n$ is cut with $q-1$ linearly independent linear forms $w_1, \dots, w_{q-1}.$ Generically, a complex $q$-dimensional variety $V^q$ intersected with such a non-degenerate set of linear forms $w_1, \dots, w_{q-1}$ will yield a collection of curves. Depending upon the position of $w_1, \dots, w_{q-1}$ and on the variety $V^q,$ the number of curves may differ and the intersection may or may not take place along the singular locus of $V^q.$ In particular, the fact that D'Angelo took an infimum over the non-degenerate sets $\{w_1, \dots, w_{q-1}\}$ means that the infimum could be realized with a set $\{w_1, \dots, w_{q-1}\}$ that exactly cuts along the singular locus of the variety $V^q$ that has the highest order of contact with the hypersurface $M.$ Catlin wished to avoid such a scenario in \cite{catlinsubell}  and defined a different notion of $q$-type $D_q$ that captured the generic behavior. Like D'Angelo he also gave a definition of $q$-type for an ideal $\I$ of holomorphic germs.

In order to understand the generic behavior, we will work in the Grassmannian $G^{n-q+1}$ of all $(n-q+1)$-dimensional complex linear subspaces of $\C^n$ at some given point $x_0.$ Each point of $G^{n-q+1}$ corresponds to a non-degenerate set of linear forms $w_1, \dots, w_{q-1},$ where the $(n-q+1)$-dimensional complex linear subspace of $\C^n$ is given by $w_1= \dots = w_{q-1}=0.$ Let $V^q$ be a germ of a $q$-dimensional complex variety at $x_0.$ Generically, for $S \in G^{n-q+1},$ $V^q \cap S$ is a finite set of irreducible holomorphic curves $V^q_{S,k}$ for $k=1, \dots, P.$ Each of them is parametrized by some open set $U_k \ni 0$ in $\C,$ so  $\gamma_S^k : U_k \rightarrow V^q_{S,k}$ with $\gamma_S^k (0)=x_0.$ D'Angelo's definitions of $\Delta_1$ and $\Delta_q$ involve computing normalized orders of contact of holomorphic curves with an element, either a holomorphic germ $f \in \oka_{x_0}$ or the defining function $r$ of a real hypersurface $M$ in $\C^n$ passing through $x_0.$ In \cite{catlinsubell} Catlin thus defined two quantities $$\tau (f, V^q \cap S) = \max_{k=1, \dots, P} \frac{ {\text ord}_0 \, {\left(\gamma^k_S\right)}^* f}{{\text ord}_0 \, \gamma^k_S}$$ and $$\tau (V^q \cap S, x_0) = \max_{k=1, \dots, P} \frac{ {\text ord}_0 \, {\left(\gamma^k_S\right)}^* r}{{\text ord}_0 \, \gamma^k_S}$$ seeking to then take a supremum over all possible germs of $q$-dimensional complex varieties at $x_0.$ For these two quantities to make sense, however, it is necessary for $V^q \cap S$ to consist just of holomorphic curves, finitely many of them and even better the same number. These were Catlin's aims. Note that germs of $q$-dimensional complex varieties at $x_0$ factor directly into Catlin's definition of $q$-type, whereas they are not directly visible in D'Angelo's definition of $q$-type. Proposition 3.1 from \cite{catlinsubell} shows that there exists an open and dense subset $\tilde W$ of $G^{n-q+1}$ such that the two quantities $\tau (f, V^q \cap S) $ and $\tau (V^q \cap S, x_0)$ stay the same for every $S \in \tilde W.$ Catlin takes out three different sets $W_1,$ $W_2,$ and $W_3$ from the Grassmannian $G^{n-q+1}$ in order to arrive at his generic set $\tilde W$ on which $P,$ the number of curves in the intersection $V^q \cap S$ is the same for all $S \in \tilde W$ and the curves $V^q_{S_a,k}$ can be smoothly parametrized via a parameter $a=(a_1, \dots, a_N),$ where $N=(n-q+1)(q-1)$ is the dimension of $\tilde W.$ Furthermore, not only is the number of intersection curves $P$ constant on $\tilde W,$ but as $S_a \in \tilde W$ varies smoothly, the intersection curves $V^q_{S_a,k}$ do as well. We shall now give more details into how Catlin found $W_1,$ $W_2,$ and $W_3.$

\subsection{Construction of $W_1$} Consider the prime ideal $\I$ in the ring $\oka_{x_0}$ of all germs of holomorphic functions that vanish on the variety $V^q.$ By the Local Parametrization Theorem, there exists a set of canonical equations for $V^q;$ see for example p.16 of \cite{gunning}. Catlin chooses a special set of coordinates such that the generators of the ideal $\I$ simultaneously satisfy the Weierstrass Preparation Theorem with respect to the variables that give the regular system of parameters $\I$ has as a prime ideal in the regular local ring $\oka_{x_0}.$ Now, the intersection $V^q \cap S$ does not behave well along the singular locus of $V^q$ and where $V^q$ does not have pure dimension $q$ since the intersection $V^q \cap S$ could consist of points rather than curves. The singular locus of $V^q$ is captured by product of the discriminants of the Weierstrass polynomials that give the canonical equations for $V^q.$ There is also a generator of $\I$ that gives the part of $V^q,$ which is not of pure dimension $q.$  Catlin constructs a conic variety $X'$ with defining equation given by the product of those discriminants and of the additional generator. $W_1$ consists of all elements of $G^{n-q+1}$ that intersect $X'.$ 

\subsection{Construction of $W_2$}  A good notion of transversality is necessary for the intersection $V^q \cap S$ to be well-behaved. $V^q \cap S$  generically yields curves on which transversality cannot be tested. To get to points, we must reduce the dimension by one, and so it makes sense to projectivize. Catlin looks at the conic variety corresponding to $V^q,$ which he calls $V'.$ This variety is defined by taking the zero locus of the homogeneous polynomials given by the leading terms of all elements of $\I.$ $V'$ still has dimension $q$ and captures the tangent cone of $V^q.$ The singularities of $V^q$ are the points where the dimension of the tangent cone jumps. Let $\tilde V$ be the projective variety in $\pj^{n-1}$ corresponding to $V',$ which we can obviously define since $V'$ is the zero locus of a collection of homogeneous polynomials. $\tilde V$ has dimension $q-1.$ To every $S \in G^{n-q+1},$ there corresponds a projective plane $\tilde S$ of dimension $n-q$ in $\pj^{n-1}.$ $\tilde V \cap \tilde S$ generically must consist of finitely many points $\tilde z^1, \dots, \tilde z^D$ with transverse intersections. That means every $\tilde z^i$ is a smooth point of $\tilde V$ and the tangent spaces satisfy the condition $T_{\tilde z^i} \tilde V \cap T_{\tilde z^i} \tilde S = 0$ for $i=1, \dots, D.$ Let $W_2$ be the subset of $G^{n-q+1}$ where this generic behavior does not happen.

\subsection{Construction of $W_3$} We used the canonical equations for $V^q$ for which the variables $z_{q+1}, \dots, z_n$ give the regular system of parameters corresponding to the pure $q$-dimensional part of the variety $V^q.$ The variable $z_q$ corresponds to the additional generator that gives the non-pure dimensional part of $V^q,$ the part that has dimension at most $q-1$ as we saw in the construction of $W_1.$ The $(n-q+1)$-dimensional linear subspace $S$ is given by the linear equations $w_i= \sum^n_{j=1} \: a^i_j z_j=0$ for $i=1, \dots, q-1,$ which have to be linearly independent. In linear algebraic terms, this condition means that there exists a $(q-1) \times (q-1)$ minor of $(a^i_j)$ of full rank. For the intersection $V^q \cap S$ to behave well, however, this $(q-1) \times (q-1)$ minor should be exactly $(a^i_j)_{1 \leq i,j \leq q-1},$ the minor with respect to the complementary variables $z_1, \dots, z_{q-1}.$ Therefore, Catlin defines $$W_3= \left\{ S \in G^{n-q+1} \:\: \Big| \:\: \det (a^i_j)_{1 \leq i,j \leq q-1}=0 \right\}.$$

\subsection{Definition of $D_q$}  The sets $W_1,$ $W_2,$ and $W_3$ are clearly closed sets in $G^{n-q+1}$ and also subvarieties of $G^{n-q+1}.$ Let $\tilde W= G^{n-q+1} \setminus (W_1 \cup W_2 \cup W_3).$ $\tilde W$ is thus an open and dense set in $G^{n-q+1}.$ We now know that for every $S \in \tilde W,$ $\tau (f, V^q \cap S)$ assumes the same value so we let $$\tau(f, V^q) = {\text gen.val}_{S \in \tilde W} \left\{ \tau (f, V^q \cap S) \right\}$$ and $$\tau(\I, V^q) = \min_{f \in \I} \tau(f, V^q).$$

\medskip
\newtheorem{catlinid}{Definition}[section]
\begin{catlinid}
\label{catlinideal}
Let $\I$ be an ideal of holomorphic germs at $x_0,$ then the Catlin $q$-type of the ideal $\I$ is given by $$D_q (\I, x_0)=\sup_{V^q}  \left\{ \tau(\I, V^q) \right\},$$ where the supremum is taken over the set of all germs of $q$-dimensional holomorphic varieties $V^q$ passing through $x_0.$
\end{catlinid}

\bigskip\noindent We also now know that $\tau (V^q \cap S, x_0)$ assumes the same value for all $S \in \tilde W,$ and so we let $$\tau(V^q, x_0) = {\text gen.val}_{S \in \tilde W} \left\{ \tau (V^q \cap S, x_0) \right\}.$$

\medskip
\newtheorem{catlinft}[catlinid]{Definition}
\begin{catlinft}
\label{catlinfinitetype}
Let $M$ be a real hypersurface in $\C^n.$ The Catlin $q$-type at $x_0 \in M$ is given by $$D_q (M, x_0)=\sup_{V^q}  \left\{ \tau (V^q, x_0) \right\},$$ where the supremum is taken over the set of all germs of $q$-dimensional holomorphic varieties $V^q$ passing through $x_0.$
\end{catlinft}

\noindent Notice that trivially $\Delta_1 (\I, x_0) = D_1 (\I , x_0)$ and $\Delta_1 (M, x_0) = D_1 (M , x_0)$ as there is only one $n$-dimensional complex plane passing through $x_0$ in $\C^n.$  Therefore, relating $\Delta_q$ and $D_q$ is only necessary for $q\geq 2.$

\section{Relationships between the D'Angelo and Catlin $q$-types}
\label{relate}

An important concept in commutative algebra that both John D'Angelo and David Catlin used in order to prove different statements about their respective notions of type is the codimension of an ideal of germs of holomorphic functions $\I$ at a point $x_0,$ $D(\I, x_0).$ Here is the exact definition:

\medskip
\begin{Definition}
For $\I$ an ideal in $\oka_{x_0},$ $$D(\I, x_0) = \dim_\C (\oka_{x_0} / \I).$$
\end{Definition}

\noindent With respect to sets of non-degenerate $q-1$ linear forms $w_1, \dots, w_{q-1}$ used in the definitions of both $\Delta_q$ and $D_q.$ the codimension of an ideal of germs of holomorphic functions behaves in a particularly good way, namely the infimum over such sets equals the generic value.

\medskip
\begin{Proposition}[\cite{bazilandreea} Proposition 2.8]
\label{multgeneric}
Let $\I$ be a proper ideal in $\oka_{x_0},$ and let $x_0 \in \C^n.$
\begin{equation*}
\inf_{\{w_1, \dots, w_{q-1} \}} D\Big( (\I, w_1, \dots, w_{q-1}), x_0 \Big)={\text gen.val}_{\{w_1, \dots, w_{q-1} \}} D\Big( (\I, w_1, \dots, w_{q-1}), x_0 \Big),
\end{equation*}
where $\{w_1, \dots, w_{q-1} \}$ is a non-degenerate set of linear forms in $\oka_{x_0},$ $(\I, w_1, \dots, w_{q-1})$ is the ideal in $\oka_{x_0}$ generated by $\I, w_1, \dots, w_{q-1},$ and the infimum and the generic value are both taken over all such non-degenerate sets $\{w_1, \dots, w_{q-1} \}$ of linear forms in $\oka_{x_0}.$ In other words, the infimum is achieved and equals the generic value.
\end{Proposition}

Vasile Brinzanescu and the author thought they had an argument that established the same for $\Delta_q$ in \cite{bazilandreea}, Proposition 2.10 for a hypersurface and Proposition 2.11 for an ideal of germs of holomorphic functions, but that turned out to be false. Martino Fassina produced an example in \cite{fassina} where they differ. Here is the Fassina example in detail:

\begin{Example}
\label{fassinaex1}
Let $\I=(z_1^3-z_3 z_2, z_2^2) \subset \oka_0^3.$ To compute $\Delta_2,$ it is necessary to take a linear form, as $q=2$ in this case, so $q-1=1.$ Let $a z_1+b z_2+ c z_3$ for $a, b, c \in \C$ be any linear form. We consider $\Delta_1 ((\I, a z_1+b z_2+ c z_3), 0).$ 

\smallskip\noindent {\bf Case 1:} Let $c \neq 0.$ Without loss of generality, we can assume $c =1.$ $\Delta_1$ is known to be invariant under biholomorphic mappings so we apply the change of variables mapping the origin to itself given by $w_1=z_1,$ $w_2=z_2,$ and $w_3=z_3+a z_1 + b z_2.$
\begin{equation*}
\begin{split}
\Delta_1 ((\I, a z_1+b z_2+  z_3), 0)&=\Delta_1((w_1^3-(w_3-a w_1 - b w_2) w_2, w_2^2, w_3),0) \\&= \Delta_1((w_1^3+a w_1 w_2 + b w_2^2, w_2^2, w_3),0) =  \Delta_1((w_1^3+a w_1 w_2  , w_2^2, w_3),0).
\end{split}
\end{equation*}
We now distinguish two cases:

\smallskip\noindent {\bf Case 1.1:} $a \neq 0.$ In this case, the supremum is achieved by the curve $\displaystyle \varphi_a (t) =\left(t ,- \frac{t^2}{a}, 0\right),$ which gives $ \Delta_1((w_1^3+a w_1 w_2  , w_2^2, w_3),0)=4= \Delta_1 ((\I, a z_1+b z_2+  z_3), 0).$

\smallskip\noindent {\bf Case 1.1:} $a = 0.$ In this case, the supremum is achieved by the curve $\varphi(t) = (t, 0, 0)$ and $ \Delta_1((w_1^3, w_2^2, w_3),0)=3= \Delta_1 ((\I, b z_2+  z_3), 0).$

\medskip\noindent {\bf Case 2:} Let $c = 0.$ In this case, $\Delta_1 ((z_1^3-z_3 z_2, z_2^2, a z_1+b z_2), 0)= + \infty$ as the curve $\varphi(t) = (0,0,t)$ has infinite order of contact with the ideal $(z_1^3-z_3 z_2, z_2^2, a z_1+b z_2).$

\smallskip\noindent Altogether, we see that  $$\Delta_2 (\I, 0) =\inf_{w_1} \Delta_1 ((\I, w_1), 0) = \inf_{a,b,c} \Delta_1 ((\I, a z_1+b z_2+ c z_3), 0)=3$$ whereas $${\text gen.val}_{w_1} \, \Delta_1 ((\I, w_1), 0) = {\text gen.val}_{a,b,c} \,\Delta_1 ((\I, a z_1+b z_2+ c z_3), 0)=4.$$

\smallskip\noindent Except for the power of $z_1$ in the first generator, this example is quite close to Example~\ref{dangeloex2}, which exhibits the maximum jump for $\Delta_1.$

\end{Example}

\begin{Example}
\label{fassinaex2}
Fassina produced a more elaborate version of the previous example so that for any $k \in \N^*$ he constructs an ideal $\I  \subset \oka_0^3$ for which the difference between the generic value and the infimum is at least $k.$ Let $m$ be such that $m^2-2m>k,$ and consider $$\I=(z_1^m-z_3 z_2, z_2^m).$$ By an argument similar to the one in the example above, Martino Fassina shows that $$\Delta_2 (\I, 0) =\inf_{w_1} \Delta_1 ((\I, w_1), 0) =m$$ whereas $${\text gen.val}_{w_1} \, \Delta_1 ((\I, w_1), 0) \geq m(m-1).$$ See Proposition 2.2. of \cite{fassina}.
\end{Example}

\noindent As a result of Fassina's examples, Vasile Brinzanescu and the author issued a correction in \cite{bazilandreeaerr}. Furthermore, it became clear that the theorems in \cite{bazilandreea} pertained to the generic notion so Vasile Brinzanescu and the author gave the following definition in \cite{bazilandreeaerr}:

\begin{Definition}
Let $2 \leq q \leq n.$ If $\I$ is an ideal in $\oka_{x_0},$
$$\tilde \Delta_q(\I, x_0)={\text gen.val}_{\{w_1, \dots, w_{q-1} \}} \:\:\ \Delta_1\Big( (\I, w_1, \dots, w_{q-1}), x_0 \Big),$$ where the generic value is taken over all non-degenerate sets $\{w_1, \dots, w_{q-1} \}$ of linear forms in $\oka_{x_0},$ $(\I, w_1, \dots, w_{q-1})$ is the ideal in $\oka_{x_0}$ generated by $\I, w_1, \dots, w_{q-1},$ and $\Delta_1$ is the D'Angelo $1$-type. Likewise, if $M$ is a real hypersurface in $\C^n$ and $x_0 \in M,$ $$\tilde \Delta_q (M, x_0) ={\text gen.val}_{\{w_1, \dots, w_{q-1} \}} \:\:\ \Delta_1\Big( (\I(M), w_1, \dots, w_{q-1}), x_0 \Big),$$ where $(\I(M), w_1, \dots, w_{q-1})$ is the ideal in $C^\infty_{x_0}$ generated by all smooth functions $\I(M)$ vanishing on $M$ along with $w_1, \dots, w_{q-1}.$ 
\end{Definition}

\begin{Remark}
Just as in Definition~\ref{qfinitetypeequiv}, the curves $\varphi \in \cv(n,x_0)$ in the definition of $ \Delta_1$ satisfy that they sit in the complex linear subspace given by $w_1 = \dots = w_{q-1}=0.$
\end{Remark}

\noindent The main results of \cite{bazilandreea} thus became:

\begin{Theorem}[Theorem 1.1 \cite{bazilandreea}, \cite{bazilandreeaerr}]
Let $\I$ be an ideal of germs of holomorphic functions at $x_0,$ then for $1 \leq q \leq n$ $$D_q(\I, x_0) \leq \tilde \Delta_q(\I, x_0) \leq \left(D_q(\I, x_0)\right)^{n-q+1}.$$
\end{Theorem}

\begin{Theorem}[Theorem 1.2 \cite{bazilandreea}, \cite{bazilandreeaerr}]
Let $\Omega$ in $\C^n$ be a domain with $\smooth$
boundary. Let $x_0 \in b \Omega$ be a point on the boundary of
the domain, and let $1 \leq q <n.$
\begin{enumerate}
\item[(i)] $D_q(b \Omega, x_0) \leq \tilde \Delta_q(b \Omega, x_0);$
\item[(ii)] If $\tilde \Delta_q(b \Omega, x_0)<\infty$ and the domain is $q$-positive at $x_0$ (the $q$ version of D'Angelo's property P), then $$\tilde \Delta_q(b \Omega, x_0) \leq 2 \left(\frac{ D_q(b \Omega, x_0)}{2} \right)^{n-q}.$$
\end{enumerate}
 In particular, if $b \Omega$ is pseudoconvex at $x_0$ and $\tilde \Delta_q(b \Omega, x_0)<\infty,$ then $$D_q(b \Omega, x_0) \leq \tilde \Delta_q(b \Omega, x_0) \leq 2 \left(\frac{ D_q(b \Omega, x_0)}{2} \right)^{n-q}.$$
\end{Theorem}

\noindent Trivially, $\Delta_1,$ $\tilde \Delta_1,$ and $D_1$ coincide and $\Delta_q \leq \tilde \Delta_q$ for $2 \leq q \leq n.$ Clearly, it was essential to effectively relate $\Delta_q$ with $\tilde \Delta_q$ and also elucidate whether the Catlin $q$-type $D_q$ and $\tilde \Delta_q$ might in fact be equal. Vasile Brinzanescu and the author carried out this work in \cite{bazilandreea2}. First of all, it follows easily from \cite{opendangelo} that $\Delta_q$ and $\tilde \Delta_q$ are effectively related for ideals of germs of holomorphic functions:

\medskip
\begin{Proposition}[Proposition 1.2 \cite{bazilandreea2}]
\label{normalandtilde}
Let $\I$ be an ideal of germs of holomorphic functions in $\oka_{x_0}.$ Let $2 \leq q \leq n,$ then $$\Delta_q(\I, x_0) \leq \tilde \Delta_q(\I, x_0) \leq \left( \Delta_q(\I, x_0) \right)^{n-q+1}.$$
\end{Proposition}

\noindent There is an effective relationship between $\Delta_q$ and $\tilde \Delta_q$ for points on the boundary of a domain as well, but it requires $q$-positivity as a hypothesis and slightly more work:

\medskip
\begin{Proposition}[Proposition 1.3 \cite{bazilandreea2}]
\label{normalandtildebdry}
Let $\Omega \subset \C^n$ be a domain with smooth boundary, and let $x_0 \in b \Omega$ be one of its boundary points. Assume $2 \leq q <n.$ If $\Omega$ is $q$-positive at $x_0,$ then
$$\Delta_q(b \Omega, x_0) \leq \tilde \Delta_q(b \Omega, x_0) \leq 2 \left( \Delta_q(b \Omega, x_0) \right)^{n-q}.$$ If the domain is pseudoconvex at $x_0,$ then it is $q$-positive, so the same inequality holds.
\end{Proposition}

\noindent Proving the equality of $D_q$ and $\tilde \Delta_q$ in \cite{bazilandreea2} was more delicate as it required a very rigorous examination of Catlin's construction of his notion of $q$-type and then defining a cylinder variety with a given curve as a base and a given complex linear subspace $W$ as its directrix subspace so that the resulting $q$-dimensional cylindrical variety would produce exactly the base curve when intersected with any element from the open and dense set $\tilde W$ in the Grassmannian $G^{n-q+1}$ from Catlin's definition of $D_q.$ The non-trivial part is that $\tilde W$ is specific to each $q$-dimensional variety as explained in Section~\ref{fincatlin}. The idea of constructing a cylinder variety was inspired by a remark made by Catlin on pp.147-8 of \cite{catlinsubell}. He observed that it would be necessary to piece together curves and create a $q$-dimensional complex variety in order to compare $\Delta_q$ and $D_q.$ Here is the lemma from \cite{bazilandreea} that shows such a cylinder variety exists:

\medskip
\begin{Lemma}[Lemma 3.1 \cite{bazilandreea2}]
\label{cylinderlemma}
Let $2 \leq q \leq n-1.$ Given a holomorphic curve $\Gamma$ passing through a point $x_0 \in \C^n$ and a $(q-1)$-dimensional hyperplane $Z$ passing through $x_0$ and satisfying that the tangent line to the curve $\Gamma$ at $x_0$ is not contained in $Z,$ there exists a germ of a $q$-dimensional cylinder $C^q$ at $x_0$ that contains the curve $\Gamma$ and whose tangent space at $x_0$ contains $Z.$
\end{Lemma}

\noindent At a singular point $x_0,$ $\Gamma$ could have have more than one tangent line. In that case, none of the tangent lines should be contained in $Z.$ This lemma while simple turns out to be a very effective tool in the study of $q$-types as will be shown in a future publication of Vasile Brinzanescu and the author. Here are now the main results of \cite{bazilandreea2}:

\begin{Theorem}[Main Theorem 1.4 \cite{bazilandreea2}]
\label{mainthm}

\begin{enumerate}[(i)]
\item Let $\I$ be an ideal of germs of holomorphic functions at $x_0,$ and let $2\leq q \leq n,$ then $$D_q(\I, x_0) = \tilde \Delta_q(\I, x_0).$$
\item Let $\Omega \subset \C^n$ be a domain with smooth boundary, let $x_0 \in b \Omega$ be one of its boundary points, and let $2 \leq q <n.$ Then $$D_q(b \Omega, x_0) = \tilde \Delta_q(b \Omega, x_0).$$
\end{enumerate}
\end{Theorem}

\smallskip
\begin{Corollary}[Corollary 1.5 \cite{bazilandreea2}]
\label{maincor}
\begin{enumerate}[(i)]
\item Let $\I$ be an ideal of germs of holomorphic functions at $x_0,$ and let $2\leq q \leq n,$ then $$\Delta_q(\I, x_0) \leq D_q(\I, x_0) \leq \left( \Delta_q(\I, x_0) \right)^{n-q+1}.$$
\item Let $\Omega \subset \C^n$ be a domain with smooth boundary, and let $x_0 \in b \Omega$ be one of its boundary points. Let $2 \leq q <n,$ and assume that $\Omega$ is $q$-positive at $x_0.$ Then $$\Delta_q(b \Omega, x_0) \leq D_q(b \Omega, x_0) \leq 2 \left( \Delta_q(b \Omega, x_0) \right)^{n-q}.$$ In particular, if the domain is pseudoconvex at $x_0,$ then the same inequality holds.
\end{enumerate}
\end{Corollary}

The consequence of the results in \cite{bazilandreea2} are twofold. First of all, the Catlin $q$-type $D_q$ becomes much easier to compute, and it is now crystal clear how it relates to the D'Angelo $q$-type. Second of all, Martino Fassina's example, Example~\ref{fassinaex1}, shows that  $$\Delta_2(\I,0)= 3 < 4 = \tilde \Delta_2(\I,0) = D_2(\I,0),$$ making it the first example of an ideal for which the D'Angelo and the Catlin $q$-types differ. It thus answers a question left open from the publication of Catlin's paper \cite{catlinsubell} in 1987. Furthermore, Fassina's Example~\ref{fassinaex2} shows that there exist ideals for which the difference between the Catlin $q$-type and the D'Angelo $q$-type is as large as we want, but in spite of that, by Corollary~\ref{maincor} the Catlin $q$-type and the D'Angelo $q$-type are still effectively related. This statement suffices for any effective computation of subelliptic gain in the $\dbar$-Neumann problem on $(p,q)$ forms. For more examples of computations of $\Delta_q$ and $D_q,$ the reader is directed to \cite{bazilandreea}.

\bibliographystyle{plain}
\bibliography{CatlinDAngeloType}

\end{document}